\documentclass[12pt,leqno]{amsart}
\usepackage{amssymb,verbatim,enumerate,ifthen,amsmath}
\usepackage{mathrsfs}
\usepackage[frak = euler]{mathalpha}
\usepackage{color}  
\usepackage{soul}  
\usepackage[mathscr]{eucal}
\usepackage[utf8]{inputenc}
\usepackage[T1]{fontenc}
\usepackage{mathrsfs}
\def\N{\mathbb{N}}
\def\R{\mathbb{R}}

\def\Z{\mathbb{Z}}
\def\T{\mathscr{T}}
\def\I{\mathscr{I}}

\newtheorem{theorem}{Theorem}[section]
\newtheorem*{theorem*}{Theorem}
\def\Thm#1#2{\ifthenelse{\equal{#1}{*}}{\begin{theorem*}#2\end{theorem*}}
             {\begin{theorem}\label{T#1}#2\end{theorem}}}
\newtheorem{Atheorem}{Theorem}

\def\thm#1{Theorem~\ref{T#1}}
\newtheorem{proposition}[theorem]{Proposition}
\newtheorem*{proposition*}{Proposition}
\def\Prp#1#2{\ifthenelse{\equal{#1}{*}}{\begin{proposition*}#2\end{proposition*}}
{\begin{proposition}\label{P#1}#2\end{proposition}}}

\newtheorem{corollary}[theorem]{Corollary}
\newtheorem*{corollary*}{Corollary}
\def\Cor#1#2{\ifthenelse{\equal{#1}{*}}{\begin{corollary*}#2\end{corollary*}}
             {\begin{corollary}\label{C#1}#2\end{corollary}}}

\newtheorem{lemma}[theorem]{Lemma}
\newtheorem*{lemma*}{Lemma}
\def\Lem#1#2{\ifthenelse{\equal{#1}{*}}{\begin{lemma*}#2\end{lemma*}}
             {\begin{lemma}\label{L#1}#2\end{lemma}}}
\def\lem#1{Lemma~\ref{L#1}}
\theoremstyle{definition}
\newtheorem{remark}[theorem]{Remark}
\newtheorem*{remark*}{Remark}
\def\Rem#1#2{\ifthenelse{\equal{#1}{*}}{\begin{remark}\rm #2\end{remark}}
             {\begin{remark}\label{R#1}\rm #2\end{remark}}}

\newtheorem{example}[theorem]{Example}
\newtheorem*{example*}{Example}
\def\Exa#1#2{\ifthenelse{\equal{#1}{*}}{\begin{example*}\rm #2\end{example*}}
             {\begin{example}\label{Ex#1}\rm #2\end{example}}}

\def\eq#1{{\rm(\ref{E#1})}}
\def\Eq#1#2{\ifthenelse{\equal{#1}{*}}
  {\begin{equation*}\begin{aligned}#2\end{aligned}\end{equation*}}
  {\begin{equation}\begin{aligned}\label{E#1}#2\end{aligned}\end{equation}}}

\begin{document}
\begin{flushright}
\end{flushright}
\vspace{5mm}

\date{\today}

\title{Decomposition of Approximately Monotone and Convex Sequences}

\author[Angshuman R. Goswami]{Angshuman R. Goswami}
\address{Department of Mathematics, University of Pannonia, 
Veszprém, Egyetem u. 10, H-8200.}
\email{goswami.angshuman.robin@mik.uni-pannon.hu}

\subjclass[2000]{Primary: . 26A51; Secondary: 39B62}
\keywords{Approximately monotone \& Convex sequences, Decomposition and approximation of sequences}
\begin{abstract}
In this paper, we primarily deal with approximately monotone and convex sequences. We start by showing that any sequence can be expressed as the difference between two nondecreasing sequences. One of these two monotone sequences act as the majorant of the original sequence, while the other possesses  non-negativity. Another result establishes that an approximately monotone(increasing) sequence can be closely approximated by a non-decreasing sequence. A similar assertion can be made for approximately convex sequence.
A sequence $\big<u_n\big>_{n=0}^{\infty}$ is said to be approximately convex (or $\varepsilon$-convex) if the following inequality holds under the mentioned assumptions
$$
u_{i}- u_{i-1}\leq u_{j}- u_{j-1}+\varepsilon\quad\quad\mbox{ where}\quad\quad i,j\in\N \quad \mbox{with} \quad i<j
.$$
We proved that an approximately convex sequence can be written as the algebraic difference of two specific types of sequences. The initial sequence contains sequential convexity property, while the other sequence possesses the Lipschitz property. Moreover, we introduce an operator $\T$, that termed as a twisting operator. In a compact interval $I(\subseteq\R)$; we characterized the convex function with this newly introduced operator.
Besides various sequence decomposition results and study related to $\T$-operator on convex function; a characterization regarding non-negative sequential convexity, a fractional inequality, implication of $\T$ operator on various types of functions, relationship between a convex function and a convex sequence are also included. \\

Motivation, previous research in this direction, various applications and linkage with some other branches of mathematics are discussed in the introduction section. 

\end{abstract}
\maketitle
\section*{Introduction} 
Functional equation and inequalities is one of the few branches of Mathematics; which is very much interdisciplinary. One of the popular research area in this field is generalizing the concept of convexity and then establishing various characterizations of it. There are hundreds of books and scientific articles from the last century that contributes in this direction. \cite{Monotonicity, Kuczma, Hyers, Becken, Ger, Gerr, Mishra, Krasniqi} and the reference therein plays an important role in the development of convex function theory and also reflects immense research possibility in this subject. However the introduction to the concept of sequential convexity is relatively new. The first mention of the term convex sequence is found in \cite{Mitrinovicc}. The importance of this terminology can be understood with the simple fact that arithmetic, Fibonacci and many other familiar sequences are only some specific types of convex sequences. Therefore, investigation of various properties of such sequences can be carried out under the umbrella of sequential convexity. Since the discovery of it, several advancements have been made and most of these are closely analogous to the ordinary convex functions. For instance; studies related to higher-order convex sequence, the establishment of a discrete version of Hermite-Hadamard type inequalities, various implications and applications of sequential convexity in the theory of trigonometric functions etc. can be found in the papers \cite{HDCS, GAB_PAL, Debnath, GauSte, Latreuch, pecaric, Essen}. Due to this reason, convex sequences are often labeled as discrete versions of convex functions.\\

A function $f:I(\subseteq\R) \to\R$ is said to be convex if for any $x,y\in I$ and for all $t\in]0,1[$; it satisfies the following functional inequality 
\Eq{*}{
f(tx+(1-t)y)\leq tf(x)+(1-t)f(y).
}
Besides this standard definition of convex function; we can define functional convexity in many other ways. The one obvious way is replacing $u:=tx+(1-t)y$ in the above inequality. After simplification, we can obtain the following equivalent form that also represents convexity of $f$ 
\Eq{2001}{
\dfrac{f(u)-f(x)}{u-x}\leq \dfrac{f(y)-f(u)}{y-u}\quad x,u,y\in I \quad \mbox{with}\quad x<u<y.
}
Both of these definitions can be used to represent/determine the convex function $f$. However, throughout this paper, we are going to extensively use \eq{2001} to establish various results.\\

The researchers have been working in the field of convex analysis for centuries. The first breakthrough in the field of approximate convexity was found in the paper \cite{Hyers}. It is shown that for any $\delta>0$, a function $f:I\to\R$, satisfying the inequality
\Eq{*}{
f(tx+(1-t)y)\leq tf(x)+(1-t)f(y)+\delta \quad  \mbox{for any} \quad x,y\in I\quad\mbox{and for all} \quad t\in[0,1]
}
can be expressed as the algebraic summation of two functions. One of these possesses convexity; while the supremum norm of the other function is less than $\delta/2.$
In the last few decades, there have been many decomposition, approximation, and characterization of various versions of approximate and higher-order convex functions. These can be found in papers like \cite{ Hardy, Hyers, Becken, Ger, Gerr, Mishra} etc. One such generalization is $\varepsilon$-convexity which was introduced by P\'ales in his paper \cite{Monotonicity}. Let $\varepsilon>0$ to be fixed. A function $f:I\to\R$ is said to be $\varepsilon$-convex if for any $x,y\in I$ and $t\in[0,1]$; the following inequality holds true
\Eq{*}{
f(tx+(1-t)y)\leq tf(x)+(1-t)f(y)+\varepsilon t(1-t)|y-x|.
} 
Later in this paper, we are going to discuss more about this special type of approximate convexity.\\

These findings open up new possibilities for studying the known definition of sequential convexity in an approximate sense and finding a linkage of it with approximate convex functions.
So far we have not found any article that closely deals with this direction. Besides, there are many scopes to provide new characterization to convex functions. In this paper, our findings are mostly related to these two topics.\\
 
In the first fold of the paper, we primarily deal with sequences. We begin with showing that any sequence $\big<u_n\big>_{n=0}^{\infty}$ can be decomposed as the sum of two sequences $\big<v_n\big>_{n=0}^{\infty}$ and $\big<w_n\big>_{n=0}^{\infty}$, where $\big<v_n\big>_{n=0}^{\infty}$ is a non-decreasing majorant of $\big<u_n\big>_{n=0}^{\infty}$, while $\big<w_n\big>_{n=0}^{\infty}$ is a non-negative and monotone(increasing) sequence.\\

Then we define the notion of an approximately convex sequence about which we have already mentioned in the abstract.
We prove that $\varepsilon$-convex sequence can be closely approximated by an ordinary convex sequence. More precisely a $\varepsilon$-convex sequence $\big<u_n\big>_{n=0}^{\infty}$ can be decomposed as the algebraic summation of a convex sequence $\big<v_n\big>_{n=0}^{\infty}$ and a controlled wavering sequence $\big<w_n\big>_{n=0}^{\infty}$  which rate of variation is bounded above by ${\varepsilon}/{2}.$  A similar type of assertion can also be made for approximately monotone sequences. We showed that a Cauchy sequence can be tightly approximated by a non-decreasing sequence in $\R$.
Moreover, we also establish that if we draw a graph by joining the points of a convex sequence with free hands; we will ultimately end up with a convex function. Mathematically, for any convex sequence $\big<u_n\big>_{n=0}^{\infty}$; there exists a continuous function $f:\R_+\to\R$; that satisfies $u_n=f(n)$  for all $n\in\N\cup\{0\}$. This more concretely justifies that the convex sequence is nothing but a discrete version of the convex function.  This simplified several results of the well-cited paper \cite{Debnath}. We will discuss it in detail later in the concerned section.
Finally, we showed that; for any arbitrarily chosen convex sequence $\big<u_n\big>_{n=0}^{\infty}$ from the domain of a non-decreasing convex function $f: I\to\R$; the corresponding functional sequence $\big<f(u_n)\big>_{n=0}^{\infty}$ is convex as well. Also, at the end of the first section, we mentioned some arising open problems regarding higher-order and approximate convex sequences.   \\

In the second part of the paper, we introduce a new operator which is termed as a twisting operator ($\T$). The implication of the twisting operator $\T$ on any function $f: [a,b]\to\R$; results a new function $\T_f$; as shown below
$$\T(f(x)):=\T_{f}(x)=f(a+b-x)\quad \mbox{for all} \quad x\in [a,b].$$
We studied changes in functional properties of monotone and periodic functions upon application of the operator $\T$. Besides,  characterizing convex functions in terms of this newly introduced operator, we provide results to some functional equations too. \\

The description of notions can be found at the very beginning of each of the sections.

\section{On monotone and convex sequences}
Throughout this section the symbols $\N$,$\,\R$ and $\R_+$ will denote the sets of natural, real and non-negative real numbers respectively.
A sequence $\big<u_n\big>_{n=0}^{\infty}$ is said to be convex if it satisfies the following discrete functional inequality 
\Eq{1}{2 u_n\leq u_{n-1}+u_{n+1} \quad \mbox{for all} \quad n\in\N.}
In other words, $\big<u_n\big>_{n=0}^{\infty}$ is said to be sequentially convex, if the sequence $\big<u_n-u_{n-1}\big>_{n=1}^{\infty}$ possesses monotonicity.
\Prp{1}{Let $\big<u_n\big>_{n=0}^{\infty}$ be a non-negative sequence such that $u_n\leq \sqrt{u_{n-1}u_{n+1}}$ holds for all $n\in\N$. Then the sequence is convex.}
\begin{proof}
Since $\big<u_n\big>_{n=0}^{\infty}$ is non-negative; we can conclude $\Big(\sqrt{u_{n+1}}-\sqrt{u_{n}}\Big)^{2}\geq 0.$ By expanding this expression and utilizing the condition on $\big<u_n\big>_{n=0}^{\infty}$, we obtain the following inequality
\Eq{*}{
2u_n\leq 2\sqrt{u_{n-1}u_{n+1}}\leq u_{n-1}+u_{n+1}.
}
This yields that $\big<u_n\big>_{n=0}^{\infty}$ is convex and completes the proof.
\end{proof}
Besides, the above proposition also shows that geometric sequences possess convexity. The next proposition gives a characterization to convex sequence under some minimal assumptions.
\Prp{2}{Let $\big<u_n\big>_{n=0}^{\infty}$ be a decreasing,  convex sequence with positive terms if and only if the reciprocal sequence $\big<v_n\big>_{n=0}^{\infty}$ defined by $v_n=\dfrac{1}{u_n}$ is positive, non-decreasing and concave.}
\begin{proof}
The non-negativity and monotonicity (increasingness)of the sequence $\big<v_n\big>_{n=0}^{\infty}$ is self-explanatory.
The convexity of the sequence $\big<u_n\big>_{n=0}^{\infty}$ yields \eq{1}. Utilizing convexity, decreasingness and non-negativity of the sequence $\big<u_n\big>_{n=0}^{\infty}$, we can obtain the following inequality
\Eq{*}{
(u_n-u_{n-1})u_{n+1}\leq (u_{n+1}-u_{n})u_{n-1}\quad\mbox{for all}\quad n\in\N.
}
Dividing both sides of the above inequality by $(u_{n-1}u_nu_{n+1})$, and rearranging the terms we get
\Eq{*}{
\dfrac{1}{u_{n-1}}+\dfrac{1}{u_{n+1}}\leq \dfrac{2}{u_n}.
}
This shows that the sequence $\big<v_n\big>_{n=0}^{\infty}$ is concave.
The reverse implication can be obtained in a similar way.
\end{proof}
The next proposition shows that any sequence can be represented as the difference of two non-deceasing sequences.
\Prp{4}{Any sequence $\big<u_i\big>_{i=0}^{\infty}$ can be expressed as the difference of two non-decreasing sequences $\big<v_i\big>_{i=0}^{\infty}$  and $\big<w_i\big>_{i=0}^{\infty}$  of $\big<u_i\big>_{i=0}^{\infty}$. Where the  sequence $\big<v_i\big>_{i=0}^{\infty}$ also act as a majorant of $\big<u_i\big>_{i=0}^{\infty}$ and the second sequence $\big<w_i\big>_{i=0}^{\infty}$ is non-negative.
}
\begin{proof}
Let $\big<u_i\big>_{i=0}^{\infty}$ be any arbitrary sequence. We can define two sequences $\big<v_i\big>_{i=0}^{\infty}$ and $\big<w_i\big>_{i=0}^{\infty}$ as follows
\Eq{*}{
v_i:=
\begin{cases}
u_0,\quad \mbox{if}\quad i=0\\
u_0+\overset{n}{\underset{i=1}\sum}|u{_i}-u_{i-1}|\quad n\in\N
\end{cases}
\quad\quad\quad\mbox{and}
\quad\quad\quad w_i=v_i-u_i
.}
We can express $u_i=v_i-w_i$ for all $i\in\N\cup\{0\}$.
By definition, It can be easily observed that 
\Eq{*}{
v_i\geq u_i \quad\mbox{and} \quad w_i\geq 0 \quad \mbox{for all} \quad  i\in\N\cup\{0\}
.}
Therefore, only remains to show the monotonicity of these two newly defined sequences.

It can be easily observed that $u_0\leq u_1$. Let $n\in \N$ to be arbitrary. Then we can compute
\Eq{*}{
v_{n+1}=u_0+\sum_{i=1}^{n+1}|u_{i}-u_{i-1}|=u_0+\sum_{i=1}^{n}|u_{i}-u_{i-1}|+|u_{n+1}-u_{n}|\geq v_n.
}
This establishes the monotonicity of the sequence $\big<v_i\big>_{i=0}^{\infty}$. We can see that
\Eq{000}{
w_1-w_0=(v_1-u_1)-(v_0-u_0)=v_1-u_1=u_0+ |u_1-u_0|-u_1\geq 0.
}
Also, for any $n\in\N$ we have the inequality 
\Eq{*}{
v_{n+1}=u_0+\sum_{i=1}^{n+1}|u_{i}-u_{i-1}|=u_0+\sum_{i=1}^{n}|u_{i}-u_{i-1}|+|u_{n+1}-u_{n}|\geq v_n+u_{n+1}-u_{n}.
}
It shows $w_{n+1}\geq w_n$ for all $n\in\N$. This together with \eq{000} proves the monotonicity of the sequence $\big<w_i\big>_{i=0}^{\infty}$ and completes the proof.
\end{proof}
A sequence  $\big<u_i\big>_{i=0}^{\infty}$ is said to be $\varepsilon$-monotone, if for any $u_m \in \big<u_i\big>_{i=0}^{\infty}$; the inequality 
\Eq{4}{u_m\leq u_n+\varepsilon \quad\mbox{for all}\quad {m<n} \quad (n\in\N).}
holds. In the next theorem, we are going to show that a $\varepsilon$-monotone sequence can be approximated by an ordinary monotone(increasing) sequence. The theorem extensively uses the decomposition result of $\varepsilon$-monotone function; which was introduced by P\'ales in his paper \cite{88}   
\Thm{3}{Let $\big<u_n\big>_{n=0}^{\infty}$ be a $\varepsilon$-monotone sequence. Then there exists a continuous, almost everywhere differentiable $\varepsilon$-monotone function $f:\R_+\to\R$ satisfying $f(n)=u_n$ for all $n\in\N\cup\{0\}.$ Additionally, there exists a monotone sequence 
$\big<v_n\big>_{n=0}^{\infty}$, such that 
\Eq{009}{
|u_n-v_n|\leq\dfrac{\varepsilon}{2}\quad \mbox{holds for all}\quad n\in\N\cup\{0\}.
}
}
\begin{proof}
To prove this theorem, We construct a function $f:\R_+\to\R$  as 
\Eq{*}{
f(x):=tu_{n-1}+(1-t)u_n\quad \mbox{;} \,\, x\in [{n-1},n]
\,\,\mbox{satisfying}\,\,x=t(n-1)+(1-t)n\quad(t\in[0,1]).}
From above, the continuity of the function $f$ is clearly visible. Besides it is evident that the constructed function is differentiable everywhere except the points at $\N\cup\{0\}$. Substituting $t=0$; at the definition of $f$, we obtain $f(n)=u_n$ for all $\N\cup\{0\}$. 
Now to show  $f$ is $\varepsilon$-monotone; let $x,y\in \R_+$ with $x<y$. It leads to the following 
\Eq{*}{x\in [{n_1-1},{n}] \quad \mbox{and} \quad y\in [{n_2-1},{n_2}],\quad \mbox{where}\quad n_1\leq n_2\quad \quad (n_1,n_2\in\N).}

Based on it $f(x)$  and $f(y)$ will satisfy the following two inequalities
\Eq{*}
{
\min\Big\{u_{n_1-1},u_{n_2}\Big\}\leq f(x)\leq \max\Big\{u_{n_1-1},u_{n_1}\Big\}\quad\mbox{and}\quad
\min\Big\{u_{n_2-1},u_{n_2}\Big\}\leq f(y)\leq \max\Big\{u_{n_2-1},u_{n_2}\Big\}
}
By, using $\varepsilon$-monotonicity; from the two inequalities above, we can compute
\Eq{*}{
f(x)-f(y)\leq \max\Big\{u_{n_1-1},u_{n_1}\Big\}- \min\Big\{u_{n_2-1},u_{n_2}\Big\}\leq \varepsilon.
}    
This yields that $f$ is $\varepsilon$-monotone and completes the first part of the theorem.

To prove the second part, we will utilize one of the results of $\varepsilon$-monotonicity that appeared in the paper \cite{Monotonicity}. The statement of the theorem can also be formulated as follows

If a function $f:\R_+\to\R$ is $\varepsilon$-monotone then there exists a monotone(increasing) function $g:\R_+\to\R$ that satisfies the inequality
\Eq{*}
{|f(x)-g(x)|\leq \dfrac{\varepsilon}{2}\quad\mbox{for all}\quad x\in I.}
Now choosing the elements of the sequence $\big<v_n\big>_{n=0}^{\infty}$ as $v_n=g(n)$ for all $n\in\N\cup\{0\}$ and using the above mentioned result, we can obtain the inequality \eq{009}. Besides, monotonicity of the function $g:\R_+\to\R$ ensures that the sequence $\big<v_n\big>_{n=0}^{\infty}$ is non-decreasing. This completes the proof and establishes the result.
\end{proof}

In the upcoming result, we will show the relationship of the Cauchy sequence with $\varepsilon$-monotonicity. A sequence $\big<u_n\big>_{n=0}^{\infty}$ is said to be Cauchy sequence if for any given $\delta\geq 0$ there exists $r\in \N$; that satisfies the following inequality under the mentioned conditions
\Eq{100}{
|u_m-u_n|\leq \delta \quad \mbox{for all} \quad m,n\geq r \quad (m,n,r\in\N).
}
\Prp{6}{Every Cauchy sequence is approximately monotone.}
\begin{proof}
To prove the proposition, we assume $\big<u_n\big>_{n=0}^{\infty}$ be a Cauchy sequence such that for a fixed $\delta\geq 0$, it satisfies \eq{100}. We consider 
\Eq{14}{\gamma =\max\Big\{\max|u_p-u_q|,\delta\Big\} \quad p,q\leq r  \quad (p,q,r\in\N).}
We take $\varepsilon=\gamma+\delta$ and will show that $\big<u_i\big>_{n=0}^{\infty}$ is $\varepsilon$-monotone. For that we have to consider three cases. At first, we assume $m,n\in\N\cup\{0\}$ with
$m<n\leq r$. Using \eq{14}, we can compute the following 
\Eq{*}{
u_m-u_n\leq \max_{p,q \leq r}|u_p-u_q|\leq\gamma<\varepsilon.
}
For $r\leq m<n,$ utilizing \eq{100} we can proceed as follows
\Eq{*}{
u_m-u_n\leq |u_m-u_n|\leq\delta\leq\varepsilon.
}
Finally for $m<r<n$, utilizing triangular inequality, we obtain the following 
\Eq{*}{
u_m-u_n\leq |u_p-u_r|+|u_r-u_n|\leq \gamma+\delta=\varepsilon.
}
The analogous inequalities under these three situations establish that $\big<u_i\big>_{n=0}^{\infty}$ possesses $\varepsilon$-monotonicity.
\end{proof}
The above result leaves us with the question, how effectively we can approximate a Cauchy sequence with a non-decreasing sequence. The next proposition shows the close relationship of a Cauchy sequence with a usual monotonic sequence.
\Cor{33}{Let $\big<u_n\big>_{n=0}^{\infty}$ be a Cauchy sequence in $\R$. Then there exists a monotone(increasing) sequence $\big<v_n\big>_{n=0}^{\infty}$ satisfying $\underset{n\to\infty}{\lim}u_n=\underset{n\to\infty}{\lim}v_n$}
\begin{proof}
Since $\big<u_n\big>_{n=0}^{\infty}$ is a Cauchy sequence in $\R$; it converges to a finite value. Without losing the generality, we assume that there exists $l\in\R$ such that $\underset{n\to\infty}{\lim}u_n=l$. Now we define a sequence $\big<u_n\big>_{n=0}^{\infty}$ as follows
\Eq{*}{
v_n=\inf_{n\leq k}v_k \qquad \qquad (n\in\N\cup\{0\})
}
By definition; it is evident that $v_n\leq u_n$ for all $n\in\N\cup\{0\}$. Also for arbitrary $n_1$ and $n_2\in\N\cup\{0\} $ with $n_1\leq n_2$; we can compute the inequality below
\Eq{*}{
v_{n_1}=\inf_{n_1\leq k}v_k\leq \inf_{n_2\leq k}v_k=v_{n_2} .
}
This yields the monotonicity of the sequence $\big<v_n\big>_{n=0}^{\infty}$. Now, we only need to show that $\underset{n\to\infty}{\lim}v_n=l$.

Since $\big<u_n\big>_{n=0}^{\infty}$ converges to $l$; it also indicates that $l$ is both limit superior and limit inferior for this sequence. By using the definition of limit inferior, we obtain the following
\Eq{*}{
\lim_{n\to\infty}v_n=\lim_{n\to\infty}\Big(\inf_{n\leq k}u_k\Big)=l\qquad\qquad (n\in\N\cup\{0\}).
}
This completes the proof.
\end{proof}
Before proceeding to the next theorem, we must establish the following lemma. 
\Lem{22}{Let $n\in\N$ be arbitrary. Then for any  $a_1,\cdots,a_n\in\R$ and $b_1,\cdots,b_n\in]0,\infty[$; the following inequality holds
\Eq{111}{
\min\bigg\{\dfrac{a_1}{b_1},\cdots, \dfrac{a_n}{b_n}\bigg\}
\leq\dfrac{a_1+\cdots+a_n}{b_1+\cdots+b_n}\leq \max\bigg\{\dfrac{a_1}{b_1},\cdots, \dfrac{a_n}{b_n}\bigg\}.}}
\begin{proof}
We will prove the theorem using the principle of mathematical induction. For $n=1$, there is nothing to prove. For $n=2$;  without loss of generality; for $a_1,a_2\in\R$ and $b_1,b_2\in]0,\infty[$, we assume that $\dfrac{a_1}{b_1}\leq \dfrac{a_2}{b_2}$. We can also express this in terms of two inequalities as follows
\Eq{*}{
a_1(b_1+b_2)\leq (a_1+a_2) b_1 \quad\mbox{and}\quad (a_1+a_2)b_2\leq (b_1+ b_2)a_2.
}
The two inequalities above together yield the following 
\Eq{67}
{\dfrac{a_1}{b_1}\leq \dfrac{a_1+a_2}{b_1+b_2}\leq \dfrac{a_2}{b_2}}
and validates the result for $n=2.$
Now we assume that the statement is true for any fixed $n\in\N$. That is \eq{111} holds true. Now, to establish the statement, we need to validate it for $n+1\in\N$. 

We assume $a_1,\cdots,a_n,a_{n+1}\in\R$ and $b_1,\cdots,b_n,b_{n+1}\in]0,\infty[$. First utilizing \eq{67} and then \eq{111}, we can compute the following inequality
\Eq{*}{
\dfrac{a_1+\cdots+a_n+a_{n+1}}{b_1+\cdots+b_n+b_{n+1}}\leq\max\bigg\{ \dfrac{a_1+\cdots+a_n}{b_1+\cdots+b_n},\dfrac{a_{n+1}}{a_{n+1}}\bigg\}\leq \max\bigg\{ \dfrac{a_1}{b_1}, \cdots ,\dfrac{a_n}{b_n}, \dfrac{a_{n+1}}{b_{n+1}} \bigg\}.}
In a similar way,   we can show the inequality below
\Eq{*}{
\min\bigg\{ \dfrac{a_1}{b_1}, \cdots ,\dfrac{a_n}{b_n}, \dfrac{a_{n+1}}{b_{n+1}} \bigg\}\leq \min \bigg\{ \dfrac{a_1+\cdots+a_n}{b_1+\cdots+b_n},\dfrac{a_{n+1}}{a_{n+1}}\bigg\}\leq \dfrac{a_1+\cdots+a_n+a_{n+1}}{b_1+\cdots+b_n+b_{n+1}}.
}
Combining the above two inequalities, we will have 
our desired result.
This completes the proof of the lemma.
\end{proof}
The origin $\Gamma$-function is directly linked to the study of factorial sequences. It states that a continuous function $\Gamma: \R_+\to\R_+$ is obtainable from the sequence $\big<n!\big>_{n=0}^{\infty}$ that possesses all the factorial properties. Similarly, the upcoming theorem shows that any discrete graph of a convex sequence can be extendable to a convex function. In other words, for every convex sequence $\big<u_n\big>_{n=0}^{\infty}$ there exists atleast one continuous convex function $f:\R_+\to\R$ that precisely approximates the sequential points. 
\Thm{15}{Let $\big<u_n\big>_{n=0}^{\infty}$ be a convex sequence. There there exists a continuous almost everywhere differentiable convex function $f:\R_{+}\to\R$ such that $f(n)=u_{n}$ for all $n\in\N\cup\{0\}.$\\
Conversely, if $f$ is convex in $\R_+$; then the sequence $<f(n)>_{n=0}^{\infty}$ is convex in $\N\cup\{0\}$.}
\begin{proof}
To prove the statement, we  define a function $f:\R_+\to\R$ as follows
\Eq{7878}{f(x):=tu_n+(1-t)u_{n+1}\quad \mbox{where} \quad x:=tn+(1-t)(n+1) \quad (t\in [0,1]\, \mbox{  and  }\, n\in \N\cup\{0\}).}
From the construction, it is clearly visible that $f$ is continuous in $\R$. Besides it is also evident that $f$ is differentiable everywhere except at the points $N\cup\{0\}$. Substituting $x=n$ in \eq{7878}, we obtain $f(n)=u_n$  for all $N\cup\{0\}$. We are only left to show the convexity of the function $f$. Before proceeding, we need to go through some tangential properties of the function $f$.

Overall, the function $f$ can be stated as the joining of line segments that are defined consecutively in between the points $n-1$ and $n$ for all $n\in\N$. For simplicity, we term these as functional line segments of $f$. 
Based on it, for any $x,y\in \R_+$ with $x<y$; we can conclude the statement as mentioned below
\Eq{200}{
&\mbox{\textbf{Slope of the functional line segment(segments) of $f$ that contains $(x,f(x))$}}\\
&\leq \mbox{\textbf{Slope of the line joining the points $(x,f(x))$ and $(y,f(y))$}}\\
&\leq \mbox{\textbf{Slope of the functional line segment(segments) of $f$ that contains $(y,f(y))$}.}
}

If both $(x,f(x))$ and $(y,f(y))$ lies in the same functional line segment, the statement is obvious.

Next we consider the case when $x\in[n-1,n[$ and $y\in]n,n+1]$, ($n\in\N$) that is $x$ and $y$; lies in two consecutive intervals.
Utilizing the sequential convexity property of $\big<u_n\big>_{n=0}^{\infty}$ (or $\big<f(n)\big>_{n=0}^{\infty}$) and basic geometry of slopes in straight lines, we can compute the inequality below
\Eq{*}{
f(n)-f(n-1)=\dfrac{f(n)-f(x)}{n-x}\leq \dfrac{f(y)-f(n)}{y-n}=f(n+1)-f(n)
.} 
The expression $\dfrac{f(y)-f(x)}{y-x}$ can also be written as $\dfrac{f(y)-f(n)+f(n)-f(x)}{(y-n)+(n-x)}$. 

Together with this, Using \eq{111} of \lem{22}; the above inequality can be extended as
\Eq{332}{
f(n)-f(n-1)\leq \dfrac{f(y)-f(x)}{y-x}\leq f(n+1)-f(n)
.} 
This yields the validity of statement \eq{200} for this particular case.
Finally we assume $x\in[n_1,n_1+1]$  and $y\in[n_2,n_2+1]$ where $n_1,n_2\in\N\cup\{0\}$ such that $n_2-n_1> 1$.
From \eq{111} of \lem{22} by utilizing \eq{332}; we can obtain the inequality below
\Eq{*}{
\dfrac{f(y)-f(x)}{y-x}&= \dfrac{\Big(f(y)-f(n_2)\Big)+\Big(f(n_2)-f(n_2-1)\Big)+\cdots +\Big(f(n_1+1)-f(x)\Big)}{\Big(y-n_2\Big)+\Big(n_2-(n_2-1)\Big)+\cdots +\Big((n_1+1)-x\Big)}\\
&\leq\max\bigg\{\dfrac{f(y)-f(n_2)}{y-n_2},\dfrac{f(n_2)-f(n_2-1)}{n_2-(n_2-1)},\cdots,\dfrac{f(n_1)-f(x)}{n_1-x}\bigg\}\\
&=\dfrac{f(y)-f(n_2)}{y-n_2}=f(n_2+1)-f(n_2).}
In a similar way, we can compute the following inequality as well
\Eq{*}{
f(n_1+1)-f(n_1)&=\dfrac{f(n_1+1)-f(x)}{(n_1+1)-x}\\
&=\min\bigg\{ \dfrac{f(n_1+1)-f(x)}{(n_1+1)-x},\dfrac{f(n_1+2)-f(n_1+1)}{(n_1+2)-(n_1+1)},\cdots,\dfrac{f(y)-f(n_2)}{y-n_2}\bigg\}\\
&\leq \dfrac{\Big( f(n_1+1)-f(x)\Big)+\Big(f(n_1+2)-f(n_1+1)\Big)+\cdots+\Big(f(y)-f(n_2)\Big)}{\Big((n_1+1)-x\Big)+\Big((n_1+2)-(n_1+1)\Big)+\cdots+\Big(y-n_2 \Big)} 
\\
&=\dfrac{f(y)-f(x)}{y-x}.}
The above two inequalities validate the statement \eq{200}.

Now we can proceed to establish the convexity of $f$. Without loss of generality, we assume $x,u,y\in\R_+$ with $x<u<y$ and $u\in[n,n+1]$ for a fixed $n\in\N\cup\{0\}$.  Then first using the last part of inequality in \eq{200} and then applying the initial part of the same inequality, we obtain the following
\Eq{*}{
\dfrac{f(u)-f(x)}{u-x}&\leq {f(n+1)-f(n)}\\
&\mbox{and}\\
f(n+1)-f(n) &\leq \dfrac{f(y)-f(u)}{y-u}.
}
Combining the above inequalities, we arrive at \eq{2001}. This shows that the function $f$ is convex and completes the first part of the proof.

To show the converse part, we assume three consecutive arbitrary $n-1,n$ and $n+1$ from $\N$. Then the convexity of $f$ yields the following inequality
\Eq{*}{
2f(n)=2f\bigg(\dfrac{(n-1)+(n+1)}{2}\bigg)\leq f(n-1)+f(n+1).
}
This establishes the sequential convexity of $\big<f_n\big>_{n=0}^{\infty}$ and completes the proof.
\end{proof}
One of the important papers of convex sequence is \cite{Debnath}, which was written by Debnath and Wu.  In Lemma 7 of that paper, without establishing the convexity of $f$, they showed that for a non-decreasing convex function $g$ on $\R$, the composite function $g\circ f:\R_+\to\R$ also possesses convexity. On the other hand, if the function $g$ is concave and decreasing then $g\circ f$ is concave on $\R_+$. Later this lemma was used to prove one of the main theorems and other subsequent results. From \thm{15} of our paper, it is clear that the function $f$ is convex. And from the well known results of composite functions, now it turns out that the lemma 7 is obvious. Some of the other results of \cite{Debnath} can also be shortened in the same manner.\\

We are now defining three approximate sequential terminologies. We will see how closely, they are interconnected with one another. 
A sequence $\big<u_n\big>_{n=0}^{\infty}$ is said to be $\varepsilon$-convex if for any $i,j\in\N\cup\{0\}$ with $i<j$; it satisfies the following inequality
\Eq{44}{u_{i}-u_{i-1}\leq u_{j}+u_{j-1}+\varepsilon.} 
It is clearly evident that a convex sequence naturally possesses $\varepsilon$-convexity. 
Substituting $j=i+1$ in \eq{44}, we can also observe that sequential $\varepsilon$-convexity leads to
the discrete inequality $2u_{i}\leq u_{i-1}+u_{i+1}+\varepsilon$ for all $i\in \N$. \\

Previously, we discussed approximately monotone sequences. A sequence $\big<v_n\big>_{n=0}^{\infty}$ is termed as $\varepsilon$-H\"older if and only if both $\big<v_n\big>_{n=0}^{\infty}$ and $\big<-v_n\big>_{n=0}^{\infty}$ are $\varepsilon$-monotone. That is the following inequality holds true
\Eq{*}{
|v_n-v_m|\leq \varepsilon \quad\mbox{for all}\quad m,n\in\N\cup\{0\}.
}

On the other hand, a sequence $\big<w_i\big>_{i=0}^{\infty}$ is said to be a Lipschitz sequence if there exists a $l\geq 0$ such that it satisfies the following inequality alongside mentioned conditions
\Eq{*}{
|w_n-w_m|\leq l|n-m| \quad\mbox{for all} \quad n,m\in \N\cup\{0\}. 
}
The infimum of all such $l$ that satisfies the above inequality is termed as Lipschitz modulo of the sequence $\big<w_i\big>_{i=0}^{\infty}.$ 
These definitions lead us to the following proposition.
\Prp{20000}{
Let $\big<u_n\big>_{n=0}^{\infty}$ and $\big<v_n\big>_{n=0}^{\infty}$ are respectively convex and ${\varepsilon}/{2}$-H\"older sequences. Then the algebraic sum of these two sequences forms a $\varepsilon$-convex sequence.
}
\begin{proof}
Let consider the sequence $\big<w_n\big>_{n=0}^{\infty}=\big<u_n+v_n\big>_{n=0}^{\infty}$ and $i,j\in\N$ such that $i<j$. The mentioned sequences give the following three inequalities
\Eq{*}{
u_i-u_{i-1}\leq u_j-u_{j-1} \qquad v_{i}-v_{i-1}\leq \varepsilon/2 \qquad \mbox{and} \qquad  - \varepsilon/2\leq v_j-v_{j-1}.
}
Summing up these inequalities side by side, we arrive at the following
\Eq{*}{
w_i-w_{i-1}\leq w_j-w_{j-1}+\varepsilon.
}
This yields that $\big<w_i\big>_{i=0}^{\infty}$ is $\varepsilon$-convex and completes the proof.
\end{proof}
The theorem below shows that a $\varepsilon$-convex sequence can be decomposed as a usual convex sequence and a Lipschitz sequence. To establish the result we will need the notion of $\varepsilon$-convex function that we have already discussed in the introduction section. As shown in \eq{2001} of the introduction; following the same substitution; an equivalent definition of functional $\varepsilon$-convexity can be formulated as follows. For a given $\varepsilon>0$; a function $f: I\to\R$ is termed as $\varepsilon$-convex if for any $x,u,y\in I$ with $x<u<y$; we have the following functional inequality
\Eq{55}{
\dfrac{f(u)-f(x)}{u-x}-\dfrac{\varepsilon}{2}\leq \dfrac{f(y)-f(u)}{y-u}+\dfrac{\varepsilon}{2}
}
Now we can proceed to the theorem which can also be seen as the Hyers-Ulam stability result in the discrete version of approximate convexity.
\Thm{9}{Let $\big<u_n\big>_{n=0}^{\infty}$ be a $\varepsilon$-convex sequence. Then there exists a $\varepsilon$-convex function $f:\R_+\to\R$ such that it satisfies $f(n)=u_n$ for all $n\in\N\cup\{0\}$. Moreover, there exists a convex sequence $\big<v_n\big>_{n=0}^{\infty}$ such that $\big<u_n-v_n\big>_{n=0}^{\infty}$ results in a Lipschitz sequence with the Lipschitz modulo bounded above by $\varepsilon$.}
\begin{proof}
To prove the statement, first we construct a function $f:\R_+\to\R$ as defined in \eq{7878} of \thm{15}. Clearly $f$ interpolates $\big<u_n\big>_{n=0}^{\infty}$ in all points of $\N\cup\{0\}$.
We will now show that $f$ is a $\varepsilon$-convex function. To establish the result, we need to show two important inequalities. Let $x,y\in\R_+$ with $x<y$ such that $x\in[n_1,n_1+1]$ and $y\in[n_2,n_2+1]$. Then
\Eq{991}{
\dfrac{f(y)-f(x)}{y-x}-\dfrac{\varepsilon}{2}&\leq \max_{n_1<n<n_2}\big\{f(n+1)-f(n)-\varepsilon/2 \big\}\quad\quad\quad(n\in\N)\\
&\mbox{and}\\
\min_{n_1<n<n_2}\big\{f(n+1)-&f(n)+\varepsilon/2 \big\}\leq \dfrac{f(y)-f(x)}{y-x}+\dfrac{\varepsilon}{2}\quad\quad\quad(n\in\N).
}
For readability purposes, we will only prove the first inequality. The establishment of the second inequality will be analogous.
With our assumptions regarding $x$ and $y$, we can rewrite the expression $\dfrac{f(y)-f(x)}{y-x}-\dfrac{\varepsilon}{2}$ as follows
\Eq{*}{
\dfrac{\Big(f(y)-f(n_2)-(y-n_2)\varepsilon/2\Big)+\cdots+\Big(f(n_1+1)-f(x)-(n_1+1-x)\varepsilon/2\Big)}{\big(y-n_2\big)+\cdots+\big(n_1+1-x\big).}
}
Using \eq{111} of \lem{22}, we can conclude the following
\Eq{*}{
\dfrac{f(y)-f(x)}{y-x}-\dfrac{\varepsilon}{2}&=\dfrac{f(y)-f(x)-(y-x)\varepsilon/2}{y-x}\\
&\leq \max\bigg\{\dfrac{f(y)-f(n_2)-(y-n_2)\varepsilon/2}{y-n_2},\cdots\\
&\quad\quad \cdots \cdots\dfrac{f(n_1+1)-f(x)-(n_1+1-x)\varepsilon/2}{n_1+1-x}
\bigg\}\\
&=\max\Big\{f(n_1+1)-f(n_1)-\varepsilon/2,\cdots\cdots,  f(n_2+1)-f(n_2)-\varepsilon/2\Big\}.}
This establishes the first part of \eq{991}. 
Now to complete the proof we assume that $x,u,y\in\R_+$ with $x<u<y$. Then according to the results obtained above; we will have distinct $i$ and $j\in\N\cup\{0\}$ with $i\leq u\leq j$  satisfying the following two inequalities
\Eq{121}{
\dfrac{f(u)-f(x)}{u-x}-\dfrac{\varepsilon}{2}&\leq f(i+1)-f(i)-\dfrac{\varepsilon}{2}\\
&\mbox{and}\\
f(j+1)-f(j)+\dfrac{\varepsilon}{2}&\leq \dfrac{f(y)-f(u)}{y-u}+\dfrac{\varepsilon}{2}.}
But $\varepsilon$-convexity of the sequence $\Big<f(n)\Big>_{n=0}^{\infty}$ yields
\Eq{*}{
f(i+1)-f(i)-\dfrac{\varepsilon}{2}\leq f(j+1)-f(j)+\dfrac{\varepsilon}{2}.
}
This together with \eq{121}, establishes \eq{55}.

Since $x,u$, and $y$ arbitrary with $x<u<y$, according to the definition of $\varepsilon$-convex function, $f$ possesses $\varepsilon$-convexity. This completes the first part of the theorem.

To prove the second part, we will utilize one of the decomposition results of functional $\varepsilon$-convexity that was formulated by Pales in his paper \cite {Monotonicity}. In view of our paper, It can be formulated as follows

"If a function $f:\R_+\to\R$ is $\varepsilon$-convex then it can be expressed as the algebraic summation of an ordinary convex function $g:\R_+\to\R$ and a Lipschitz function $h:\R_+\to\R$  such that Lipschitz modulo $lip(h)$ is bounded 
above by $\varepsilon.$ That is 
$$\underset{x,y\,\in\,\R_+}{\sup}\left|\dfrac{f(y)-f(x)}{y-x}\right|\leq \varepsilon.$$"
Now we choose the elements of the sequence $\big<v_n\big>_{n=0}^{\infty}$ as $v_n=g(n)$ for all $n\in\N\cup\{0\}$. From the converse part of \thm{9999}, we can conclude that $\big<v_i\big>_{n=0}^{\infty}$ is a convex sequence. Using the stated result above, we can obtain the inequality \eq{009} as well. This proves all the assertions in the theorem and completes the proof.
\end{proof}
\Thm{9999}{Let $f:\R_+\to\R$ be a non-decreasing convex function.  Then for any strictly non-decreasing convex sequence $\big<u_n\big>_{n=0}^{\infty} \in \R_+$; the corresponding functional sequence $\big<f(u_n)\big>_{n=0}^{\infty}$ is convex. }
\begin{proof}
Since $f$ is convex; then for any non-decreasing convex sequence $\big<u_n\big>_{n=0}^{\infty} \in \R_+$; it will satisfy the following inequality
\Eq{*}{
 \dfrac{f(u_n)-f(u_{n-1})}{u_n-u_{n-1}}\leq \dfrac{f(u_{n+1})-f(u_{n})}{u_{n+1}-u_{n}}\qquad\mbox{for all} \quad n\in\N.
}
Again using the monotonicity convexity property of the convex sequence $\big<u_n\big>_{n=0}^{\infty}$, we obtain
\Eq{*}{
u_{n}-u_{n-1}\leq u_{n+1}-u_n \qquad\mbox{for all} \quad n\in\N.
}
Multiplying the above two inequalities side by side and utilizing the non-negativity of all the involved terms, we arrive at the following 
\Eq{*}{
f(u_{n})-f(u_{n-1})\leq f(u_{n+1})-f(u_n) \qquad\mbox{for all} \quad n\in\N.
}
This yields the monotonicity of the sequence $\big<f(u_{n})-f(u_{n-1}\big>_{n=0}^{\infty}$; in other words establishes that $\big<f(u_n)\big>_{n=0}^{\infty}$ is convex.
\end{proof}
Higher-order convexity plays an important role in the development of function theory. To define the notion of higher-order convexity, we first need to become familiar with the concept of divided difference (see the papers \cite{Becken, Ger, Gerr, Mishra} and the referred articles in those). In the inductive way, the first and the $n^{th}$ ordered divided difference of the function $f:I(\subseteq\R)\to\R$ can be represented as follows
\Eq{*}{
[u_i,u_{i+1};f]&=\dfrac{f(u_i)-f(u_{i+1})}{u_i-u_{i+1}} \qquad \quad  {i\in\{0,1,\cdots,n\}}\qquad \quad u_i\in I \\
[u_0,\cdots,u_{n};f]&=\dfrac{[u_0\cdots,u_{n-1};f]-[u_1\cdots,u_{n};f]}{u_{n}-u_0} \qquad \qquad u_0,\cdots,u_{n}\in I.
}
A function $f:I\to\R$ can be termed as $n$-convex if the following assertion satisfies 
$$[u_0,\cdots,u_{n};f]\geq 0 \qquad \qquad  u_0,\cdots,u_{n}\in I.
$$
From the definition, It can be easily observed that the $1$-convex denotes the ordinary monotonicity(increasingness) of the function $f$. While the $2^{nd}$-order convexity is our well-known ordinary functional convexity.

A sequence $\big<u_i\big>_{i=0}^{\infty}$ is said to be $n^{th}$-order convex, if it satisfies the below mentioned non-negativity
\Eq{*}{
\sum_{r=0}^{n}(-1)^r\binom{n}{r}u_{n-r}\geq 0 \qquad \qquad 0\leq r\leq n\quad (r,n\in\N).
}
Studies related to higher-ordered sequential convexity can be found in the papers \cite
It is evident that in sequence the first order convexity denotes the usual monotonicity. And one can also see that the $2^{nd}$-order convexity is our usual sequential convexity. 

It is clearly visible that  $n^{th}$ order sequential convexity is just the discrete version of $n^{th}$-order convex function. Now the very basic question that we can formulate is, "For a given $n^{th}$-order convex sequence, is it possible to define a continuous convex function $f:\R_+\to\R$ of the same order that passes through the sequence? Another important open question is the possibility of approximating a $\varepsilon$- $n^{th}$ order convex sequence by an ideal $n^{th}$ order convex sequence. \\

In the next section, we are going to study an operator and its various implications in monotone, periodic and convex functions.
\section{The Twisting Operator}
Throughout this section, $I$ will denote the non-empty and non-singleton interval $[a,b]$.
For a function $f: I\to\R$; the operator $\T$ will be termed as twisting operator and defined as bellow
$$\T({f}(x)):=f(a+b-x)\quad \mbox{for all}\quad x\in I.$$ For the upcoming results $\T_f$ will denote the newly generated function after implication of the twisting operator $\T$ on $f$.

It can be easily observed that under the respective functional operations, the set consists of the identity operator $\I$; defined as $\I(x):=x$ (for all $x\in I$) and the twisting $\T$ together forms a group. More specifically $\{\I,\T\}$ is an commutative group. One can also observe the following 
\Eq{6543}{
\T^2(f)=\T(\T(f))=T_{\T_{f}}=\I \quad \mbox{that is}\quad \T^2(f)(x)=x \quad\quad (x\in I).
}
The first proposition below shows the effect of the twisting operator in monotone functions.
\Prp{7}{Let $f:I\to\R$ be a non-decreasing function. Then $\T_f$ will generate a decreasing function}
\begin{proof}
To prove the statement, we assume $x,y\in I$ with $x<y$. This implies the inequality
$a+b-y\leq a+b-x$. Now using the monotonicity(increasingness) of $f$, we can obtain
\Eq{*}{
\T_f(y)=f(a+b-y)\leq f(a+b-x)=\T_f(x).
} 
The above inequality establishes decreasingness of $\T_f$ and completes the proof.
\end{proof}
A function $f:I\to\R$ is said to be periodic with a period $L$, where $L<b-a$; if for any $x\in I$, there exists some $n\in\Z$ such that the inequality $f(x)=f(x+nL)$ holds.

The stated proposition below shows the impact twisting operator on periodic function.
\Prp{55}{ Let $0<L<\ell(I)$, where $\ell(I)$ denotes the length of the interval $I$.
If $f:I\to\R$ is a periodic function with the period $L$, then $\T_f$ also results in a periodic function.
}
\begin{proof}
To prove this statement, we assume $x\in I$ be arbitrary. Then due to our assumption regarding $\ell(I)$; there must exists some $n\in\Z$ such that $x+nL\in I.$ This implies for all those selective $n$, the elements $a+b-(x+nL)'s$ and $a+b-x \in I$.
Using the $L$-periodicity of $f$, we can compute the following
\Eq{*}{
\T_f(x)=f(a+b-x)=f(a+b-(x+nL))=\T_f(x+nL). 
}
Since $x$ is arbitrary, we can conclude that $\T_f$ is periodic. This completes the proof.
\end{proof}
Next, we will study two functional equations and the underlying results in them. Before proceeding to these propositions, we must mention that $I'$ will denote the compact interval  $\bigg[-\dfrac{b-a}{2}, \dfrac{b-a}{2}\bigg].$ This first proposition shows that if both $f$ and $\T_f$ represent the same function, then it is almost equivalent to the evenness property of the function. 
\Prp{8}{Let $f:I\to\R$ satisfies the functional equation $f(x)=f(a+b-x)$. Then the function $g:I'\to\R$ defined by $g(\theta):=f\bigg(\dfrac{a+b}{2}+\theta\bigg)$ is even.}
\begin{proof}
To prove the proposition, we assume $\theta\in I'$ to be arbitrary. This implies both $\bigg(\dfrac{a+b}{2}-\theta\bigg)$ and $\left(\dfrac{a+b}{2}+\theta\right)\in I$. By using the substitution, $x=\dfrac{a+b}{2}-\theta$ in the functional equation $f(x)=f(a+b-x)$, we obtain the following 
$$f\bigg(\dfrac{a+b}{2}-\theta\bigg)=f\bigg(\dfrac{a+b}{2}+\theta\bigg).$$
This shows that $g(-\theta)=g(\theta)$ and validates the statement.
\end{proof}
The proposition below deals with the scenario when $f+T_f$=0.
\Prp{*}{Let $f:I\to\R$ satisfies the functional equation $f(x)+f(a+b-x)=0$. Then the function $g:I'\to\R$ defined by $g(\theta):=f\bigg(\dfrac{a+b}{2}+\theta\bigg)$ is even.}
\begin{proof}
The proof of the proposition is analogous to the previous result. Hence omitted.
\end{proof}
In the upcoming theorem, we provide a characterization of the convex function in terms of the twisting operator $\T$. we can also see that the twisting operator results in a convex function when it acts on a convex function.

\Thm{11}{Let $f:I\to\R$. Then the following conditions are equivalent to each other:
\begin{enumerate}[(i)]
\item $f$ is convex;
\item For all $x,u,y\in I$ with $x<u<y$, 
 \Eq{*}{
 \dfrac{(\T_f)(u)-(\T_f)(x)}{u-x}\leq \dfrac{(\T_f)(y)-(\T_f)(u)}{y-u};
 }
\item There exists a function $\varphi: I\to\R$ such that for all $x,u\in I$,
 \Eq{ux}{
 (\T_f)(u)+(x-u)\varphi(u)\leq (\T_f)(x);
 }
\item For all $n\in\N$, $x_1,\dots,x_n\in I,$ $t_1,\dots,t_n\geq0$ with $ t_1+\dots+t_n=1$,
 \Eq{*}{
 (\T_f)(t_1x_1+\dots+t_nx_n)\leq t_1\big(\T_f)(x_1)+\cdots+t_n(\T_f)(x_n).
 }
 \item  $\T_f$ is convex
\end{enumerate}
}

\begin{proof} {\it(i)$\Rightarrow$(ii):} Assume that $f$ is $\Phi$ convex and let $x<u<y$ to be arbitrary elements of $I$. This implies $a+b-y<a+b-u<a+b-x$. We choose $t\in[0,1]$ such that $a+b-u=t(a+b-y)+(1-t)(a+b-x)$. In other words, we select $\T=\dfrac{u-x}{y-x}$; and hence \eq{1} can be rewritten as follows
\Eq{*}{
  f(a+b-u)\leq\dfrac{u-x}{y-x}f(a+b-y)+\dfrac{y-u}{y-x}f(a+b-x).
}
Therefore, 
\Eq{*}{
  (y-u+u-x)(\T_f)(u)\leq&(u-x)(\T_f)(y)+(y-u)(\T_f)(x)\\
}
Rearranging the above inequality we obtain, 
\Eq{*}{
 \dfrac{(\T_f)(u)-(\T_f)(x)}{u-x}&\leq \frac{(\T_f)(y)-(\T_f)(u)}{y-u}.
}\\

{\it(ii)$\Rightarrow$(iii):}
Assume that {\it(ii)} holds and define the function $\varphi$ on $I$ by
\Eq{*}{
  \varphi(u):=\underset{x<u}{\mathrm{\sup}}\bigg(\dfrac{(\T_f)(u)-(\T_f)(x)}{u-x}\bigg) \qquad(u\in I).
}
In view of condition {\it(ii)}, for all $x<u<y$ in $I$, we have
\Eq{6}{
\dfrac{(\T_f)(u)-(\T_f)(x)}{u-x}\leq\varphi(u)\leq \dfrac{(\T_f)(y)-(\T_f)(u)}{y-u}.
}
From the left hand side inequality in \eq{6}, we get
\Eq{*}{
(\T_f)(u)+(x-u)\varphi(u)\leq (\T_f)(x) \qquad (x\in I,\,x<u).
}
Similarly, from the right hand side inequality in \eq{6} (replacing $y$ by $x$, it follows that
\Eq{*}{
(\T_f)(u)+(x-u)\varphi(u)\leq (\T_f)(x) \qquad (x\in I,\,u<x).
}
Now, combining the two inequalities above, the condition {\it(iii)} follows (also in the case $x=u$).\\

{\it(iii)$\Rightarrow$(iv):} 
To deduce {\it(iv)} from {\it(iii)}, let $x_1,\dots,x_n\in I, $ and $t_1,\dots,t_n\geq 0$ with $t_1+\dots+t_n=1$ and $u:=t_1x_1+\dots+t_nx_n$.  Substituting $x$ by $x_i$ in the inequality of condition {\it(iii)} we will obtain $n$ inequalities. Then multiplying this inequality by $t_i$, and finally adding up the inequalities side by side, we get
\Eq{*}{
\sum_{i=1}^{n} t_i\big((\T_f)(u)+(x_i-u)\varphi(u)\big)\leq&\sum_{i=1}^{n} t_i\big((\T_f)(x_i).
}
Using that $\sum_{i=1}^{n} t_i(x_i-u)=0$, the above inequality simplifies to the inequality of condition {\it(iv)}.\\

{\it(iv)$\Rightarrow$(v):} To deduce the $\Phi$-convexity of $f$ from condition {\it(iv)}

let $x,y\in I$ and $t\in[0,1]$. Taking $n=2$, $x_1:=x$, $x_2:=y$, $t_1:=t$ and $t_2:=1-t$ in condition {\it(iv)}, it is immediate to see that the inequality reduces to 
\Eq{*}{
(\T_f)(tx+(1-t)y)\leq t(\T_f)(x)+(1-t)(\T_f)(y).
}
This yields the convexity of the function $\T_f$.\\

{\it(v)$\Rightarrow$(vi):}
The following the same steps from $\it(i)$ to $\it(v)$; and using \eq{6543} we can  finally obtain $\T_{f}$, that is $f$ is convex.\\

This establishes equivalency among all the implications and completes the proof.
\end{proof}
Motivated by the condition {\it(iii)} of \thm{1}, we say that $\varphi:I\to\R$ is a \emph{$\Phi$-slope function for $\T_f$} if it satisfies inequality \eq{ux} for all $x,u\in I$. The next proposition establishes the monotonic property of $\phi$.

\Prp{14}{The slope function $\varphi$ defined in \eq{ux} is nondecreasing.}
\begin{proof}
To prove the statement, we assume that $x,u\in I$ with $x<u$. 
Then interchanging the roles of $x$ and $u$ in the \eq{ux}(third assertion of the above theorem), we get the following
\Eq{*}{
(\T_f)(x)+(u-x)\varphi(x)\leq (\T_f)(u).
} 
Now adding up the above inequality with \eq{ux} side by side, we arrive at
\Eq{*}{(u-x)(\varphi(x)-\varphi(u))\leq 0.}
This implies $\Phi(u)\leq \Phi(x)$ and due to the arbitrariness of $u$ and $x$; the monotonicity is established.
\end{proof}
\bibliographystyle{plain}

\end{document}